\documentclass[10pt]{amsart}
\usepackage{amsfonts,amssymb,amscd,amsmath,enumerate,verbatim,calc}
\usepackage{amsthm}
\usepackage{url}
\newtheorem{thrm}{Theorem}[section]
\newtheorem{remark}[thrm]{Remark}
\newtheorem{Example}[thrm]{Example}
\numberwithin{equation}{section}

\author{Ramesh  Kumar Muthumalai}
\address{Ramesh  Kumar Muthumalai\newline
 Department of Mathematics\newline
 D.G.Vaishnav College,\newline
 Arumbakkam, Chennai-106,\newline
  Tamil Nadu, India.}
\email{ramjan\_80@yahoo.com. }
\begin{document}
\title[Generalization of Ramanujan Method ]{Generalization of Ramanujan Method of Approximating  root of an equation}
\begin{abstract}
We generalize Ramanujan method of approximating the smallest root of an equation which is found in Ramanujan Note books, Part-I. We provide simple analytical proof to study convergence of this method. Moreover, we study iterative approach of this method on approximating a root with arbitrary order of convergence.
\end{abstract}
\maketitle
\section{Introduction.}
Ramanujan method of approximating the smallest root $z_0$ of an equation of the form, 
\begin{equation}
\sum_{k=1}^\infty  A_kz^k=1
\tag{1.1}
\end{equation}
is found in Chapter 2 of Ramanujan note books, part-I  without proof.
It is assumed that all other roots of (1.1) have moduli strictly greater than 
$|z_0|$. For $z$ sufficiently small, write
\begin{equation}
\frac{1}{1-\sum_{k=1}^\infty  A_kz^k}=\sum_{k=1}^\infty P_kz^{k-1}
\tag{1.2}
\end{equation}
It follows easily that $P_1=1$ and 
\begin{equation}
P_n=\sum_{j=1}^{n-1}A_jP_{n-j}, \hspace{1cm} n\geq 1 
\tag{1.3}
\end{equation}
Ramanujan gives (1.3) and claims, with no hypothesis, that $P_n/P_{n+1}$ approaches a root of (1.1)\cite{1,5}. In this present study, we  give analytic proof of his method and generalize this to  approximate a root of nonlinear equation with the arbitrary order of convergence. 

\section{Generalization of Ramanujan method.} Suppose that  $f$  be a  real-valued function, defined  and continuous on a bounded closed interval $[a,b]$ of the real line that satisfies  $f(a)f(b)<0$. Assume, further, that  $f$ is  continuously differentiable  on  $(a,b)$. We wish to find real number $\alpha\in[a,b]$ such that $f(\alpha)=0$. Let  $z$ be an approximate of $\alpha$ that satisfies $|z-\alpha|<1$ and $f'(z)\neq 0$. Define  
\begin{align}
P(z)f(z)=1,\hspace{1cm} \mbox{If }z\neq\alpha
\tag{2.1}
\end{align}
Now applying Leibniz's rule\cite{3} of $n^{th}$ derivative on (2.1) and after simplification, we find that
\begin{align}
P_n(z)=-\frac{1}{f(z)}\sum_{k=0}^{n-1}\binom{n}{k}P_k(z)f^{(n-k)}(z),\hspace{1cm}n>1
\tag{2.2}
\end{align}
where $P_k(z)$ is $k^{th}$ derivative of $1/f(z)$, $(k=1,2,\ldots,n)$ and $P_0(z)=1/f(z)$. Using (2.2), the next three values of $P(z)$'s are listed as follow as 
\begin{equation*}
P_1(z)=-\frac{f'(z)}{f(z)^2}\qquad\qquad\qquad\qquad\quad
\end{equation*}
\begin{equation*}
P_2(z)=-\frac{f^{(2)}(z)}{f(z)^2}+2\frac{f'(z)^2}{f(z)^3}
\qquad\qquad
\end{equation*}
\begin{equation*}
\qquad\quad P_3(z)=-\frac{f^{(3)}(z)}{f(z)^2}+6\frac
{f'(z)f^{(2)}(z)}{f(z)^3}-6\frac{f'(z)^3}{f(z)^4}
\end{equation*}
and so on. Now define for $n\geq 1$
\begin{equation*}
H_n(z)=\frac{D^{n}\left(\frac{z-\alpha}{f(z)}\right)}{D^{n}\left(\frac{1}{f(z)}\right)}
\tag{2.3}
\end{equation*}
The operator $D^{n}$ in (2.3) denotes  $n^{th}$ order differentiation.
Assume that,  $\alpha$ is a zero of  order one for $f(z)$ then $f(\alpha)=0$ and  $f'(\alpha)\neq 0$.
Set
\begin{equation*}
Q(z)=\frac{1}{\sum_{j=1}^\infty\frac{f^{(j)}(\alpha)}{j!}(z-\alpha)^{j-1}}
\tag{2.4}
\end{equation*} 
It  follows easily from  Taylor's series \cite{3} that  $Q(z)=\frac{z-\alpha}{f(z)}$ for $z\neq\alpha$ and  $Q(\alpha)=1/f'(\alpha)$. Since $\alpha$ is zero of order one of $f(z)$, $Q(z)$ is differentiable at $\alpha$. 
Subtituting (2.4) in (2.3), we find that
\begin{align*}
H_n(z)&=\frac{D^{n}Q(z)}{D^{n}\left(\frac{Q(z)}{z-\alpha}\right)}
\qquad\qquad\qquad
\end{align*}
Expanding $Q(z)$ by Taylor's series  at $z=\alpha$ on both numerator and denominator and applying the operator $D^n$  and after simplification, we have
\begin{align*}
\qquad\qquad\qquad\qquad&=\frac{(z-\alpha)^{n+1}\sum_{j=0}^\infty\frac{Q^{(n+j)}
(\alpha)}{j!}(z-\alpha)^{j}}{n!(-1)^nQ(\alpha)+\sum_{j=1}^\infty \frac{Q^{(n+j)}(\alpha)}{n+j}\frac{(z-\alpha)^{n+j}}{(j-1)!}}
\end{align*}
Thus, we obtain
\begin{align}
H_n(z)&=O\left((z-\alpha)^{n+1}\right)\qquad\qquad\qquad
\tag{2.5}
\end{align}
Since $|z-\alpha|<1$ and letting $n \rightarrow \infty$, we obtain
\begin{align}
\lim_{n\rightarrow 0}H_n(z)=0\qquad\qquad\qquad\qquad
\tag{2.6}
\end{align}
On the other hand, using Leibnitz rule on (2.3), we find that 
\begin{equation*}
H_n(z)=\frac{(z-\alpha)D^{n}\left(\frac{1}{f(z)}\right)+nD^{n-1}\left(\frac{1}
{f(z)}\right)}{D^{n}\left(\frac{1}{f(z)}\right)}
\end{equation*}
Since $D^n(1/f(z))=P_n(z)$ and after simplfiiction, we obtain
\begin{align}
H_n(z)&=z-\alpha+n\frac{P_{n-1}(z)}{P_{n}(z)}
\tag{2.7}
\end{align}
Letting $n\rightarrow \infty$ and using (2.7), we obtain
\begin{equation}
\alpha=z+\lim_{n \rightarrow \infty} n\frac{P_{n-1}(z)}{P_{n}(z)}
\tag{2.8}
\end{equation}
\begin{remark} If we replace  $P_n(z)$ by $n!P_n$, $P_{n-1}(z)$ by $(n-1)!P_{n-1}$ and assuming $z=0$  and $f(0)=-1$  in (2.1), (2.2) and (2.8), we get Ramanujan method of approximating the smallest root of an equation of the form (1.1) \cite{3}. 
\end{remark}
The recursion formula given in (2.2) may lead numerical instability in machine computation due to divisibility of $f(z)$ when $f(z)\rightarrow 0$. To avoid such instability, multiply  $f(z)^n$ on both sides of (2.1), (2.2) and  set $T_k(z)=P_k(z)f(z)^{k+1},$ we obtain $T_0(z)=1$ and
\begin{align}
T_n(z)=-\sum_{k=0}^{n-1}\binom{n}{k}T_k(z)f(z)^{n-k-1}f^{(n-k)}(z),\hspace{0.5cm}n>1
\tag{2.9}
\end{align}
Also, the equation (2.8) becomes 

\begin{equation}
\alpha=z+\lim_{n\rightarrow \infty} nf(z)\frac{T_{n-1}(z)}{T_{n}(z)}
,\hspace{0.5cm}n=1,2,\ldots
\tag{2.10}
\end{equation}
Thus, $n^{th}$ convergent of (2.10) gives
\begin{equation*}
 \alpha\approx z+nf(z)\frac{T_{n-1}(z)}{T_{n}(z)},\hspace{0.5cm}n=1,2,\ldots
\end{equation*}
\section{Some rational approximations.}Ramanujan gives six examples   to illustrate performance of his method to approximate the smallest root of equations.  The approximations are given  in the rational form.  In particular, he provides an interesting rational approximation that $\log 2=375/541$. In this section, we present such rational approximations to $m^{th}$ root of some integers  and logarithm of some rational numbers.
\subsection{  $\mathbf{m^{th}}$ root of a rational number.}
Let us consider  $f(z)=z^m-a$ and $c$ be an integer which is very near to $m^{th}$ root of $a$.  Now, using (2.2), we obtain the following recursion
\begin{align}
T_n(c)&=-\sum_{k=0}^{r}\binom{n}{k}\frac{m!}{(m-n+k)!}T_k(c)\left(c^m-a\right)^{n-k-1}c^{m-n+k}
\tag{3.1}
\end{align}
and
\begin{align}
\sqrt[m]{a}=c+n(c^m-a)\frac{T_{n-1}(c)}{T_{n}(c)}
\tag{3.2}
\end{align}
where $r=n-1\mbox{ if $n\leq m$}$ and  $r=m\quad\mbox{if $n>m$}$.
The sequence of rational approximations of $\sqrt{2}$ are 
\begin{equation*}
\frac{99}{70}, \frac{239}{169},  \frac{577}{408}, \frac{1393}{985}, \frac{3363}{2378}, \frac{8119}{5741},\frac{19601}{13860},  \frac{47321}{33461}                   
\end{equation*}
which is corected up to  3, 4, 4, 5, 7, 7, 8, 9  decimal places respectively. Table 1 lists   rational approximations of some irrational numbers using MATLAB.
\begin{table}[h]
\centering
\caption{Rational approximations of Irrational numbers.}
\begin{tabular}{c|c|c|c}
\hline
S.N0 & Irrational No. & Rational Approx. & significant digits\\
\hline
1&$\sqrt[3]{9}$&$\frac{50623}{24337}$&10-digits\\
2&$\sqrt[9]{511}$&$\frac{4603}{2302}$ & 9-digits\\
3&$\sqrt[3]{2}$&$\frac{6064}{4813}$&8-digits\\
4 & $\sqrt[5]{3100}$&$\frac{3110}{623}$ & 7-digits\\
\hline
\end{tabular}
\label{tab:1}
\end{table}
\subsection{Some Logarithmic values.}Let us consider  $f(z)=\log(z+1)-a$ and $z=0$ be an initial approximation (where $a$ is a rational number). Table 2 lists rational approximations of logarithmic  numbers using MATLAB.
\begin{table}[ht]
\centering
\caption{Rational approximations of logarithmic values.}
\begin{tabular}{c|c|c|c}
\hline
S.No & Logarithm Nos. & Rational Approx. & significance digits\\
\hline
1&$\log_e 1.5$&$\frac{3858}{9515}$ & 7-digits\\
2&$\log_e 2.0$&$\frac{32781}{47293}$ & 6-digits\\
3&$\log_e 3.0$&$\frac{12667}{11530}$&7-digits\\
4& $\log_e 1.2$&$\frac{724}{3971}$ & 6-digits\\
\hline
\end{tabular}
\label{tab:2}
\end{table}
\subsection{Some other examples.} Ramanujan lists the first ten convergents to the real root of $x+x^3=1$, with the last convergent being $13/19=0.684210\ldots$. By Newton's method, the root is $0.682327804\ldots$.   Here, we provide three examples to approximate  root of  $x^3-2x-5=0$,  $e^x-3=0$ and $x=\sin x +\frac{1}{2}$ by taking the initial approximations $x_0=2$, $x_0=1$ and $x_0=1$ respectively. Table 3 lists  first ten convergents of roots all the three equations. The roots of first two equations converge to 14 decimal places at $9^{th}$ and $8^{th}$ convergent respectively whereas the root of $x=\sin x +\frac{1}{2}$ converges
slowly  to  $1.49730038909589$ at $23^{th}$ convergent. We see  from these examples, this method of approximating a real root is  very slow process. However, we get  good approximation when $n$ is large.
\begin{table}[ht]
\centering
\footnotesize{
\caption{Examples for Generalized Ramanujan method}
\begin{tabular}{c|c|c|c}
\hline
It & $x^3-2x-5$=0 & $e^x-3=0$& $x=\sin x +\frac{1}{2}$\\
\hline
0 &  2.00000000000000 &  1.00000000000000 &  1.00000000000000\\
1 &  2.10000000000000 &  1.10363832351433 &  1.58288042035629\\
2 &  2.09433962264151 &  1.09853245432531 &  1.51838510578857\\
3 &  2.09455842997324 &  1.09861223692174 &  1.50077867834371\\
4 &  2.09455128205128 &  1.09861230157476 &  1.49783013943789\\
5 &  2.09455148653822 &  1.09861228868606 &  1.49735888023541\\
6 &  2.09455148143875 &  1.09861228866513 &  1.49730334991792\\
7 &  2.09455148154375 &  1.09861228866810 &  1.49729959647640\\
8 &  2.09455148154234 &  1.09861228866811 &  1.49730005778495\\
9 &  2.09455148154232 &  1.09861228866811 &  1.49730030987454\\
\hline
\end{tabular}}
\label{tab:3}
\end{table}
\section{ Iterartive approach on Generalized Ramanujan method.}
In this section, we  study iterative approach on  generalized Ramanujan's method  to improve the accuracy of the root of nonlinear equations. If $z_0$ be an initial approximate root of $f(z)$ then using (2.10), 
\begin{equation}
z_{m+1}=z_m+ nf(z_m)\frac{T_{n-1}(z_m)}{T_{n}(z_m)}\hspace{1 cm}m=0,1,2,\ldots
\tag{4.1}
\end{equation}
where
\begin{align}
T_n(z_m)=-\sum_{k=0}^{n-1}\binom{n}{k}T_k(z_m)f(z_m)^{n-k-1}f^{(n-k)}(z_m),
\hspace{0.3cm}n>1
\tag{4.2}
\end{align}
Set $z_{m+1}-\alpha=\epsilon_{m+1}$ and $z_{m}-\alpha=\epsilon_{m}$. Then by using  (2.5) and (2.7), we find that the  order of convergence of (4.1) is $n+1$. (i.e)
\begin{equation}
\epsilon_{m+1}=O\left(\epsilon_m^{n+1}\right)
\tag{4.3}
\end{equation}
Also, the condition of convergence is
\begin{equation}
\left|n+1-n\frac{T_{n-1}(z)T_{n+1}(z)}{T_{n}(z)^2}\right|< 1
\tag{4.4}
\end{equation}
Setting $n=1$ in (4.1) and after simplification, we obtain Newton-Raphson method
\begin{equation*}
z_{m+1}=z_m-\frac{f(z_m)}{f'(z_m)},\hspace{1 cm}m=0,1,2,\ldots
\end{equation*}
Setting $n=2$ in (4.1) and after simplification, we obtain Halley's method
\begin{equation*}
z_{m+1}=z_m-\frac{f(z_m)/f'(z_m)}{1-\frac{1}{2}f(z_m)f''(z_m)/f'(z_m)^2}
\hspace{1 cm}m=0,1,2,\ldots
\end{equation*}
By varying the values of $n$, we can find more formulas with the higher order of accuracy. Thus, this new iterative  method is more generalization of Newton-Raphson method \cite{3} and Halley's method\cite{4}. 
\begin{Example} The equation $f(z)=z-\cos z$ has  exactly one  root $\alpha=0.73908513321516$ between $0$ and $1$ and starting with $z_0=0$, generate the sequence 
$z_1, z_2,\ldots$ by taking $m=1:5$ and $n=1:4$ in  (4.5)    are listed in the Table 4.
\end{Example}
\begin{table}[h]
\centering
\footnotesize{
\caption{Computation of root of $z-\cos z=0$}
\begin{tabular}{c|c|c|c|c}
\hline
$m$&$n=1$&$n=2$&$n=3$&$n=4$\\
\hline
 1&1.00000000000000  & 0.66666666666667  & 0.75000000000000  &  0.73846153846154\\
 2& 0.75036386784024 &  0.73903926244631 &  0.73908513352403 &  0.73908513321516\\
 3& 0.73911289091136 &  0.73908513321515 &  0.73908513321516 &  \\
 4& 0.73908513338528 &  0.73908513321516 &  &  \\
 5& 0.73908513321516 &   & &  \\
\hline
\end{tabular}}
\label{tab:4}
\end{table}
We observed that from Table 4, when $m=1$ and $n=1:4$ (i.e. first row) the sequence corresponds to generalized Ramaujan method studied in section 2, which converges very slowly towards the root.  Simliarly, when $n=1$ and $m=1:5$ (i.e first column) the sequence corresponds to Newton-Raphson method, which converges to root at $5^{th}$ iteration, second column  correspond to Halley method, which converges to root at $4^{th}$ iteration  and so on. 
\section{ Conclusion.} In conclusion we note that the  generalization of  Ramanujan method of approximating the real root of nonlinear equations has been developed in this article.   Firstly, we have generalized Ramanujan's method of approximating the smallest root of  equations, which is found in Ramanujan Note books - I, to any real root with simple analytic proof.  Secondly, we  have proved  that  iterative approch of this method has arbitrary order of convergence.  Moreover, we have shown that Newton-Raphson and Halley's method are special cases of this generalized Ramanujan  method. 

\end{document}